\documentclass[12pt]{amsart}
\usepackage{amssymb}

\usepackage{t1enc}
\usepackage[latin2]{inputenc}
\usepackage{verbatim}
\usepackage{amsmath,amsfonts,amssymb,amsthm}
\usepackage[mathcal]{eucal}
\usepackage{enumerate}
\usepackage[centertags]{amsmath}
\usepackage{graphics}

\newtheorem{theorem}{Theorem}
\newtheorem{lemma}{Lemma}

\newtheorem{corollary}{Corollary}

\numberwithin{equation}{subsection}

\begin{document}
\author{George Tephnadze}
\title[logarithmic means ]{The maximal operators of logarithmic means of one-dimensional
Vilenkin-Fourier series. }
\address{G. Tephnadze, Department of Mathematics, Faculty of Exact and Natural
Sciences, Tbilisi State University, Chavchavadze str. 1, Tbilisi 0128,
Georgia}
\email{giorgitephnadze@gmail.com}
\date{}
\maketitle

\begin{abstract}
The main aim of this paper is to investigate $\left( H_{p},L_{p}\right) $%
-type inequalities for maximal operators of logarithmic means of
one-dimensional Vilenkin-Fourier series.
\end{abstract}

\textbf{2000 Mathematics Subject Classification.} 42C10.

\textbf{Key words and phrases:} Vilenkin system, Logarithmic means,
martingale Hardy space.

\section{INTRODUCTION}

In one-dimensional case the weak type inequality

\begin{equation*}
\mu \left( \sigma ^{*}f>\lambda \right) \leq \frac{c}{\lambda }\left\|
f\right\| _{1}\text{ \qquad }\left( \lambda >0\right)
\end{equation*}
can be found in Zygmund \cite{Zy} for the trigonometric series, in Schipp
\cite{Sc} for Walsh series and in Pál, Simon \cite{PS} for bounded Vilenkin
series. Again in one-dimensional, Fujji \cite{Fu} and Simon \cite{Si2}
verified that $\sigma ^{*}$ is bounded from $H_{1}$ to $L_{1}$. Weisz \cite
{We2} generalized this result and proved the boundedness of $\sigma ^{*}$
from the martingale space $H_{p}$ to the space $L_{p}$ for $p>1/2$. Simon
\cite{Si1} gave a counterexample, which shows that boundedness does not hold
for $0<p<1/2.$ The counterexample for $p=1/2$ due to Goginava ( \cite{Go},
see also \cite{BGG2}).

Riesz` s logarithmic means with respect to the trigonometric system was
studied by a lot of autors.We mentioned, for instance, the paper by Szasz
\cite{Sz} and Yabuta \cite{Ya}. this means with respect to the Walsh and
Vilenkin systems by Simon\cite{Si1} and Gát\cite{Ga1}.

Móricz and Siddiqi\cite{Mor} investigates the approximation properties of
some special Nörlund means of Walsh-Fourier series of $L_{p}$ function in
norm. The case when $q_{k}=1/k$ is excluded, since the methods of Móricz and
Siddiqi are not applicable to Nörlund logarithmic means. In \cite{Ga2} Gát
and Goginava proved some convergence and divergence properties of the Nö%
rlund logarithmic means of functions in the class of continuous functions
and in the lebesque space $L_{1}.$ Among there, they gave a negative answer
to the question of Móricz and Siddiqi\cite{Mor}. Gát and Goginava\cite{Ga3}
proved that for each measurable function $\phi \left( u\right) =\circ \left(
u\sqrt{\log u}\right) $there exists an integrable function $f,$ such that

\begin{equation*}
\underset{G_{m}}{\int }\phi \left( \left| f\left( x\right) \right| \right)
d\mu \,\left( x\right) <\infty
\end{equation*}
and there exist a set with positive measure, such that the Walsh-logarithmic
means of the function diverge on this set.

The main aim of this paper is to investigate $(H_{p},L_{p})$-type
inequalities for the maximal operators of Riesz and Nörlund logarithmic
means of one-dimensional Vilenkin-Fourier series. We prove that the maximal
operator $R^{*}$ is bounded from the Hardy space $H_{p}$ to the space $L_{p}$
when $p>1/2.$We also shows that when $0<p\leq 1/2$ there exists a martingale
$f\in H_{p},$ for which

\begin{equation*}
\left\| R^{*}f\right\| _{L_{p}}=+\infty .
\end{equation*}

For the Nörlund logarithmic means we prove that when $0<p\leq 1$ there
exists a martingale $f\in H_{p}$ for which

\begin{equation*}
\left\| L^{*}f\right\| _{L_{p}}=+\infty .
\end{equation*}

Analogical theorems for Walsh-Paley system is proved in \cite{goginava}.

\section{DEFINITIONS AND NOTATIONS}

Let $N_{+}$ denote the set of the positive integers, $N:=N_{+}\cup \{0\}.$
Let $m:=(m_{0,}m_{1....})$ denote a sequence of the positive integers not
less than 2. Denote by $Z_{m_{k}}:=\{0,1,...m_{k}-1\}$ the addition group of
integers modulo $m_{k}$.

Define the group $G_{m}$ as the complete direct product of the groups $%
Z_{m_{i}}$ with the product of the discrete topologies of $Z_{m_{j}}`$ s.

The direct product $\mu $ of the measures

\begin{equation*}
\mu _{k}\left( \{j\}\right) :=1/m_{k},\text{ \qquad }(j\in Z_{m_{k}})
\end{equation*}
is the Haar measure on $G_{m_{k\text{ }}},$ with $\mu \left( G_{m}\right)
=1. $

If $\sup\limits_{n}m_{n}<\infty $, then we call $G_{m}$ a bounded Vilenkin
group. If the generating sequence $m$ is not bounded then $G_{m}$ is said to
be an unbounded Vilenkin group. \textbf{In this paper we discuss bounded
Vilenkin groups only.}

The elements of $G_{m}$ represented by sequences

\begin{equation*}
x:=\left( x_{0},x_{1},...,x_{j},...\right) ,\text{ }\left( x_{i}\in
Z_{m_{j}}\right) .
\end{equation*}

It is easy to give a base for the neighborhood of $G_{m}$

\begin{eqnarray*}
I_{0}\left( x\right) &:&=G_{m}, \\
I_{n}(x) &:&=\{y\in G_{m}\mid y_{0}=x_{0},...y_{n-1}=x_{n-1}\},\text{ }%
\,\,\left( x\in G_{m},n\in N\right) .
\end{eqnarray*}

Denote $I_{n}:=I_{n}\left( 0\right) ,$ for $n\in N_{+}.$

If we define the so-called generalized number system based on $m$ in the
following way :

\begin{equation*}
M_{0}:=1,\text{ }M_{k+1}:=m_{k}M_{k},\,\,\,(k\in N),
\end{equation*}
then every $n\in N,$ can be uniquely expressed as $n=\sum_{j=0}^{\infty
}n_{j}M_{j},$ where $n_{j}\in Z_{m_{j}},$ $(j\in N_{+})$ and only a finite
number of $n_{j}`$ s differ from zero.

Next, we introduce on $G_{m}$ an ortonormal system which is called the
Vilenkin system. At first define the complex valued function $r_{k}\left(
x\right) :G_{m}\rightarrow C,$ The generalized Rademacher functions as

\begin{equation*}
r_{k}\left( x\right) :=\exp \left( 2\pi ix_{k}/m_{k}\right) ,\text{ }\left(
i^{2}=-1,x\in G_{m},\text{ }k\in N\right) .
\end{equation*}

Now define the Vilenkin system$\,\,\,\psi :=(\psi _{n}:n\in N)$ on $G_{m}$
as:

\begin{equation*}
\psi _{n}(x):=\prod\limits_{k=0}^{\infty }r_{k}^{n_{k}}\left( x\right)
,\,\,\,\left( n\in N\right) .
\end{equation*}
Specifically, we call this system the Walsh-Paley one if $m\equiv 2.$

The Vilenkin system is orthonormal and complete in $L_{2}\left( G_{m}\right)
.$\cite{AVD,Vi}

Now we introduce analogues of the usual definitions in Fourier-analysis. If $%
f\in L_{1}\left( G_{m}\right) $ we can establish the Fourier coefficients,
the partial sums of the Fourier series, the Fej\'er means, the Dirichlet
kernels with respect to the Vilenkin system $\psi $ in the usual manner:

$\qquad $%
\begin{eqnarray*}
\widehat{f}\left( k\right) &:&=\int_{G_{m}}f\overline{\psi }_{k}d\mu ,\,\,%
\text{\qquad }\left( k\in N\right) , \\
S_{n}f &:&=\sum_{k=0}^{n-1}\widehat{f}\left( k\right) \psi _{k},\text{
\qquad }\left( n\in N_{+},S_{0}f:=0\right) , \\
\sigma _{n}f &:&=\frac{1}{n}\sum_{k=0}^{n-1}S_{k}f\,,\,\text{\qquad }\left(
n\in N_{+}\right) , \\
D_{n} &:&=\sum_{k=0}^{n-1}\psi _{k}\,,\,\,\qquad \left( n\in N_{+}\right) .
\end{eqnarray*}

Recall that
\begin{equation*}
D_{M_{n}}\left( x\right) =\left\{
\begin{array}{l}
M_{n},\text{ }\,\,\text{if\thinspace \thinspace \thinspace }x\in I_{n}, \\
0,\text{ }\,\,\text{if}\,\,x\notin I_{n}.
\end{array}
\right.
\end{equation*}

The norm (or quasinorm) of the space $L_{p}(G_{m})$ is defined by \qquad

\begin{equation*}
\left\| f\right\| _{p}:=\left( \int_{G_{m}}\left| f(x)\right| ^{p}d\mu
(x)\right) ^{\frac{1}{p}},\qquad \left( 0<p<\infty \right) .
\end{equation*}

The $\sigma -$algebra generated by the intervals $\left\{ I_{n}\left(
x\right) :x\in G_{m}\right\} $ will be denoted by $\digamma _{n}\left( n\in
N\right) .$ Denote by $f=\left( f^{\left( n\right) },n\in N\right) $ a
martingale with respect to $\digamma _{n}\left( n\in N\right) .$(for details
see e.g. \cite{We1}).

The maximal function of a martingale $f$ is defined by

\begin{equation*}
f^{*}=\sup_{n\in N}\left| f^{(n)}\right| .
\end{equation*}

In case $f\in L_{1}\left( G_{m}\right) ,$ the maximal functions are also be
given by

\begin{equation*}
f^{*}\left( x\right) =\sup\limits_{n\in N}\frac{1}{\mu \left( I_{n}\left(
x\right) \right) }\left| \int\limits_{I_{n}\left( x\right) }f\left( u\right)
d\mu \left( u\right) \right| .
\end{equation*}

For $0<p<\infty $ the Hardy martingale spaces $H_{p}$ $\left( G_{m}\right) $
consist of all martingale for which

\begin{equation*}
\left\| f\right\| _{H_{p}}:=\left\| f^{*}\right\| _{L_{p}}<\infty .
\end{equation*}

If $f\in L_{1}\left( G_{m}\right) ,$ then it is easy to show that the
sequence $\left( S_{M_{n}}\left( f\right) :n\in N\right) $ is a martingale.

If $f=\left( f^{\left( n\right) },\text{ }n\in N\right) $ is martingale then
the Vilenkin-Fourier coefficients must be defined in a slightly different
manner:

\begin{equation*}
\widehat{f}\left( i\right) :=\lim_{k\rightarrow \infty
}\int_{G_{m}}f^{\left( k\right) }\left( x\right) \overline{\Psi }_{i}\left(
x\right) d\mu \left( x\right) .
\end{equation*}

The Vilenkin-Fourier coefficients of $f\in L_{1}\left( G_{m}\right) $ are
the same as those of the martingale $\left( S_{M_{n}}\left( f\right) :n\in
N\right) $ obtained from $f$.

In the literature, there is the notion of Riesz` s logarithmic means of the
Fourier series. The n-th Riesz$^{,}$ s logarithmic means of the Fourier
series of an integrable function $f$ is defined by

\begin{equation*}
R_{n}f\left( x\right) :=\frac{1}{l_{n}}\overset{n}{\underset{k=1}{\sum }}%
\frac{S_{k}f\left( x\right) }{k},
\end{equation*}
where

\begin{equation*}
l_{n}:=\overset{n}{\underset{k=1}{\sum }}\left( 1/k\right) .
\end{equation*}

Let \{$q_{k}:$ $k>0$\} be a sequence of nonnegative numbers. The n-th Nö%
rlund means for the Fourier series of $f$ is defined by

\begin{equation*}
\frac{1}{Q_{n}}\overset{n}{\underset{k=1}{\sum }}q_{n-k}S_{k}f,
\end{equation*}
where

\begin{equation*}
Q_{n}:=\overset{n}{\underset{k=1}{\sum }}q_{k}.
\end{equation*}

If $q_{k}=k,$ then we get Nörlund logarithmic means

\begin{equation*}
L_{n}f\left( x\right) :=\frac{1}{l_{n}}\overset{n}{\underset{k=1}{\sum }}%
\frac{S_{k}f\left( x\right) }{n-k}.
\end{equation*}
It is a kind of $^{,,}$reverse$^{,,}$ Riesz` s logarithmic means.

In this paper we call this means logarithmic means.

For the martingale $f$ we consider the following maximal operators of

\begin{eqnarray*}
R^{*}f\left( x\right) &:&=\underset{n\in N}{\sup }\left| R_{n}f\left(
x\right) \right| , \\
L^{*}f\left( x\right) &:&=\underset{n\in N}{\sup }\left| L_{n}f\left(
x\right) \right| , \\
\sigma ^{*}f\left( x\right) &:&=\sup_{n\in N}\left| \sigma _{n}f(x)\right| .
\end{eqnarray*}

A bounded measurable function $a$ is p-atom, if there exists a dyadic
interval I, such that

\begin{equation*}
\left\{
\begin{array}{l}
a)\qquad \int_{I}ad\mu =0, \\
b)\ \qquad \left\| a\right\| _{\infty }\leq \mu \left( I\right) ^{-1/p}, \\
c)\qquad \text{supp}\left( a\right) \subset I.\qquad
\end{array}
\right.
\end{equation*}

\section{FORMULATION OF MAIN RESULT}

\begin{theorem}
Let $p>1/2.$ Then the maximal operator $R^{*}$ is bounded from the Hardy
space $H_{p}$ to the space $L_{p}$.
\end{theorem}

\begin{theorem}
Let $0<p\leq 1/2$. Then there exists a martingale $f\in H_{p}$ such that
\begin{equation*}
\left\| R^{*}f\right\| _{p}=+\infty .
\end{equation*}
\end{theorem}

\begin{corollary}
Let $0<p\leq 1/2$. Then there exists a martingale $f\in H_{p}$ such that
\end{corollary}

\begin{equation*}
\left\| \sigma ^{*}f\right\| _{p}=+\infty .
\end{equation*}

\begin{theorem}
Let $0<p\leq 1$. Then there exists a martingale $f\in L_{p}$ such that
\end{theorem}

\begin{equation*}
\left\| L^{*}f\right\| _{p}=+\infty .
\end{equation*}

\section{AUXILIARY PROPOSITIONS}

\begin{lemma}
\cite{We3} A martingale $f=\left( f^{\left( n\right) },n\in N\right) $ is in
$H_{p}\left( 0<p\leq 1\right) $ if and only if there exist a sequence $%
\left( a_{k},k\in N\right) $ of p-atoms and a sequence $\left( \mu _{k},k\in
N\right) $ of a real numbers such that for every n$\in N:$

\begin{equation}
\underset{k=0}{\overset{\infty }{\sum }}\mu _{k}S_{M_{n}}a_{k}=f^{\left(
n\right) },  \label{1}
\end{equation}
\end{lemma}

\begin{equation*}
\qquad \sum_{k=0}^{\infty }\left| \mu _{k}\right| ^{p}<\infty .
\end{equation*}
Moreover,

\begin{equation*}
\left\| f\right\| _{H_{p}}\backsim \inf \left( \sum_{K=0}^{\infty }\left|
\mu _{k}\right| ^{p}\right) ^{1/p},
\end{equation*}
where the infimum is taken over all decomposition of $f$ of the form (\ref{1}%
).

\section{PROOF OF THE THEOREM}

\textbf{Proof of theorem 1.} Using Abel transformation we obtain

\begin{equation*}
R_{n}f\left( x\right) =\frac{1}{l_{n}}\overset{n-1}{\underset{j=1}{\sum }}%
\frac{\sigma _{j}f\left( x\right) }{j+1}+\frac{\sigma _{n}f\left( x\right) }{%
l_{n}},
\end{equation*}
Consequently,

\begin{equation}
L^{*}f\leq c\sigma ^{*}f.  \label{2}
\end{equation}

On the other hand Weisz\cite{We2} proved that $\sigma ^{*}$ is bounded from
the Hardy space $H_{p}$ to the space $L_{p}$ when $p>1/2.$ Hence, from (2)
we conclude that $R^{*}$ is bounded from the martingale Hardy space $H_{p}$
to the space $L_{p}$ when $p>1/2.$

\textbf{Proof of theorem 2. }Let $\left\{ \alpha _{k}:k\in N\right\} $ be an
increasing sequence of the positive integers such that

\qquad
\begin{equation}
\sum_{k=0}^{\infty }\alpha _{k}^{-p/2}<\infty ,  \label{3}
\end{equation}

\begin{equation}
\sum_{\eta =0}^{k-1}\frac{\left( M_{2\alpha _{\eta }}\right) ^{1/p}}{\sqrt{%
\alpha _{\eta }}}<\frac{\left( M_{2\alpha _{k}}\right) ^{1/p}}{\sqrt{\alpha
_{k}}},  \label{4}
\end{equation}

\begin{equation}
\frac{\left( M_{2\alpha _{k-1}}\right) ^{1/p}}{\sqrt{\alpha _{k-1}}}<\frac{%
M_{\alpha _{k}}}{\alpha _{k}^{3/2}}.  \label{5}
\end{equation}

We note that such an increasing sequence $\left\{ \alpha _{k}:k\in N\right\}
$ which satisfies conditions (\ref{3})-(\ref{5}) can be constructed.

Let \qquad
\begin{equation*}
f^{\left( A\right) }\left( x\right) =\sum_{\left\{ k;\text{ }2\alpha
_{k}<A\right\} }\lambda _{k}a_{k},
\end{equation*}
where
\begin{equation*}
\lambda _{k}=\frac{m_{2\alpha _{k}}}{\sqrt{\alpha _{k}}}
\end{equation*}
and

\begin{equation*}
a_{k}\left( x\right) =\frac{M_{2\alpha _{k}}^{1/p-1}}{m_{2\alpha _{k}}}%
\left( D_{M_{2\alpha _{k}+1}}\left( x\right) -D_{M_{_{2\alpha _{k}}}}\left(
x\right) \right) .
\end{equation*}

It is easy to show that

\begin{eqnarray*}
\left\| a_{k}\right\| _{\infty } &\leq &\frac{M_{2\alpha _{k}}^{1/p-1}}{%
m_{2\alpha _{k}}}M_{2\alpha _{k}+1} \\
&\leq &(M_{2\alpha _{k}})^{1/p}=(\text{supp}\left( a_{k}\right) )^{-1/p},
\end{eqnarray*}

\begin{equation}
S_{M_{A}}a_{k}\left( x\right) =\left\{
\begin{array}{l}
a_{k}\left( x\right) ,\text{ }2\alpha _{k}<A, \\
0,\text{ }2\alpha _{k}\geq A.
\end{array}
\right.  \label{6}
\end{equation}
\begin{eqnarray*}
f^{\left( A\right) }\left( x\right) &=&\sum_{\left\{ k;\text{ }2\alpha
_{k}<A\right\} }\lambda _{k}a_{k}=\underset{k=0}{\overset{\infty }{\sum }}%
\lambda _{k}S_{M_{A}}a_{k}\left( x\right) , \\
\text{supp}(a_{k}) &=&I_{2\alpha _{k}}, \\
\int_{I_{2\alpha _{k}}}a_{k}d\mu &=&0.
\end{eqnarray*}

from (\ref{3}) and lemma 1 we conclude that $f=\left( f^{\left( n\right)
},n\in N\right) \in H_{p}.$

Let

\begin{equation*}
q_{A}^{s}=M_{2A}+M_{2s}-1,\qquad A>S.
\end{equation*}

Then we can write

\begin{eqnarray}
R_{q_{\alpha _{k}}^{s}}f\left( x\right) &=&\frac{1}{l_{q_{\alpha _{k}}^{s}}}%
\underset{j=1}{\overset{q_{\alpha _{k}}^{s}}{\sum }}\frac{S_{j}f\left(
x\right) }{j}  \label{7} \\
&=&\frac{1}{l_{q_{\alpha _{k}}^{s}}}\underset{j=1}{\overset{M_{2\alpha
_{k}}-1}{\sum }}\frac{S_{j}f\left( x\right) }{j}  \notag \\
&&+\frac{1}{l_{q_{\alpha _{k}}^{s}}}\underset{j=M_{2\alpha _{k}}}{\overset{%
q_{\alpha _{k}}^{s}}{\sum }}\frac{S_{j}f\left( x\right) }{j}  \notag \\
&=&I+II.  \notag
\end{eqnarray}

It is easy to show that

\begin{equation}
\widehat{f}(j)=\left\{
\begin{array}{l}
\frac{M_{2\alpha _{k}}^{1/p-1}}{\sqrt{\alpha _{k}}},\,\,\text{ if \thinspace
\thinspace }j\in \left\{ M_{2\alpha _{k}},...,\text{ ~}M_{2\alpha
_{k}+1}-1\right\} ,\text{ }k=0,1,2..., \\
0,\text{ \thinspace \thinspace \thinspace if \thinspace \thinspace
\thinspace }j\notin \bigcup\limits_{k=1}^{\infty }\left\{ M_{2\alpha
_{k}},...,\text{ ~}M_{2\alpha _{k}+1}-1\right\} \text{ .}
\end{array}
\right.  \label{8}
\end{equation}

Let $j<M_{2\alpha _{k}}.$Then from (\ref{4}) and (\ref{8}) we have

\begin{eqnarray}
&&\left| S_{j}f\left( x\right) \right|  \label{9} \\
&\leq &\sum_{\eta =0}^{k-1}\sum_{v=M_{2\alpha _{\eta }}}^{M_{2\alpha _{\eta
}+1}-1}\left| \widehat{f}(v)\right|  \notag \\
&\leq &\sum_{\eta =0}^{k-1}\sum_{v=M_{2\alpha _{\eta }}}^{M_{2\alpha _{\eta
}+1}-1}\frac{M_{2\alpha _{\eta }}^{1/p-1}}{\sqrt{\alpha _{\eta }}}  \notag \\
&\leq &c\sum_{\eta =0}^{k-1}\frac{M_{2\alpha _{\eta }}^{1/p}}{\sqrt{\alpha
_{\eta }}}\leq \frac{cM_{2\alpha _{k-1}}^{1/p}}{\sqrt{\alpha _{k-1}}}.
\notag
\end{eqnarray}

Consequently

\begin{eqnarray}
\left| I\right| &\leq &\frac{1}{l_{q_{\alpha _{k}}^{s}}}\underset{j=1}{%
\overset{M_{2\alpha _{k}}-1}{\sum }}\frac{\left| S_{j}f\left( x\right)
\right| }{j}  \label{10} \\
&\leq &\frac{c}{\alpha _{k}}\frac{M_{2\alpha _{k-1}}^{1/p}}{\sqrt{\alpha
_{k-1}}}\sum_{j=1}^{M_{2\alpha _{k}}-1}\frac{1}{j}  \notag \\
&&c\frac{M_{2\alpha _{k-1}}^{1/p}}{\sqrt{\alpha _{k-1}}}.  \notag
\end{eqnarray}

Let $M_{2\alpha _{k}}\leq j\leq q_{\alpha _{k}}^{s}.$ Then we have the
following

\begin{eqnarray}
S_{j}f\left( x\right) &=&\sum_{\eta =0}^{k-1}\sum_{v=M_{2\alpha _{\eta
}}}^{M_{2\alpha _{\eta }+1}-1}\widehat{f}(v)\psi _{v}\left( x\right)
+\sum_{v=M_{2\alpha _{k}}}^{j-1}\widehat{f}(v)\psi _{v}\left( x\right)
\label{11} \\
&=&\sum_{\eta =0}^{k-1}\frac{M_{2\alpha _{\eta }}^{1/p-1}}{\sqrt{\alpha
_{\eta }}}\left( D_{M_{_{2\alpha _{\eta }+1}}}\left( x\right)
-D_{M_{_{2\alpha _{\eta }}}}\left( x\right) \right)  \notag \\
&&+\frac{M_{2\alpha _{k}}^{1/p-1}}{\sqrt{\alpha _{k}}}\left( D_{j}\left(
x\right) -D_{M_{_{2\alpha _{k}}}}\left( x\right) \right) .  \notag
\end{eqnarray}

This gives that

\begin{eqnarray}
II &=&\frac{1}{l_{q_{\alpha _{k}}^{s}}}\underset{j=M_{2\alpha _{k}}}{%
\overset{q_{\alpha _{k}}^{s}}{\sum }}\ \frac{1}{j}\left( \sum_{\eta =0}^{k-1}%
\frac{M_{2\alpha _{\eta }}^{1/p-1}}{\sqrt{\alpha _{\eta }}}\left(
D_{M_{_{2\alpha _{\eta }+1}}}\left( x\right) -D_{M_{_{2\alpha _{\eta
}}}}\left( x\right) \right) \right)  \label{12} \\
&&+\frac{1}{l_{q_{\alpha _{k}}^{s}}}\frac{M_{2\alpha _{k}}^{1/p-1}}{\sqrt{%
\alpha _{k}}}\sum_{j=M_{2\alpha _{k}}}^{q_{\alpha _{k}}^{s}}\frac{\left(
D_{_{j}}\left( x\right) -D_{M_{_{2\alpha _{k}}}}\left( x\right) \right) }{j}
\notag \\
&=&II_{1}+II_{2}.  \notag
\end{eqnarray}

To discuss $II_{1}$, we use (\ref{4}). Thus we can write:

$.$

\begin{equation}
\left| II_{1}\right| \leq c\sum_{\eta =0}^{k-1}\frac{M_{2\alpha _{\eta
}}^{1/p}}{\sqrt{\alpha _{\eta }}}\leq \frac{cM_{2\alpha _{k-1}}}{\sqrt{%
\alpha _{k-1}}}.  \label{13}
\end{equation}

Since

\begin{equation}
D_{j+M_{2\alpha _{k}}}\left( x\right) =D_{M_{2\alpha _{k}}}\left( x\right)
+\psi _{_{M_{2\alpha _{k}}}}\left( x\right) D_{j}\left( x\right) ,\text{
\qquad when \thinspace \thinspace }j<M_{2\alpha _{k}},  \label{14}
\end{equation}
for $II_{2}$ we have

\begin{eqnarray}
II_{2} &=&\frac{1}{l_{q_{\alpha _{k}}^{s}}}\frac{M_{2\alpha _{k}}^{1/p-1}}{%
\sqrt{\alpha _{k}}}\sum_{j=0}^{M_{2s}}\frac{D_{j+M_{2\alpha _{k}}}\left(
x\right) -D_{M_{2\alpha _{k}}}\left( x\right) }{_{j+M_{2\alpha _{k}}}}
\label{15} \\
&=&\frac{1}{l_{q_{\alpha _{k}}^{s}}}\frac{M_{2\alpha _{k}}^{1/p-1}}{\sqrt{%
\alpha _{k}}}\psi _{_{M_{2\alpha _{k}}}}\sum_{j=0}^{M_{2s}-1}\frac{%
D_{j}\left( x\right) }{_{j+M_{2\alpha _{k}}}}.  \notag
\end{eqnarray}

We write

\begin{equation*}
R_{q_{\alpha _{k}}^{s}}f\left( x\right) =I+II_{1}+II_{2},
\end{equation*}
Then by (\ref{5} ), (\ref{7}), (\ref{10}) and (\ref{12})-(\ref{15}) we have

\begin{eqnarray*}
\left| R_{q_{\alpha _{k}}^{s}}f\left( x\right) \right| &\geq &\left|
II_{2}\right| -\left| I\right| -\left| II_{1}\right| \\
&\geq &\left| II_{2}\right| -c\frac{M_{\alpha _{k}}}{\alpha _{k}^{3/2}} \\
&\geq &\frac{c}{\alpha _{k}}\frac{M_{2\alpha _{k}}^{1/p-1}}{\sqrt{\alpha _{k}%
}}\left| \sum_{j=0}^{M_{2s}-1}\frac{D_{j}\left( x\right) }{_{j+M_{2\alpha
_{k}}}}\right| -c\frac{M_{\alpha _{k}}}{\alpha _{k}^{3/2}}.
\end{eqnarray*}

Let $0<p\leq 1/2,$ $x\in I_{2s}\backslash I_{2s+1}$ for $s=\left[ 2\alpha
_{k}/3\right] ,...,\alpha _{k}.$ Then it is evident

\begin{equation*}
\left| \sum_{j=0}^{M_{2s}-1}\frac{D_{_{j}}\left( x\right) }{_{j+M_{2\alpha
_{k}}}}\right| \geq \frac{cM_{2s}^{2}}{M_{2\alpha _{k}}}.
\end{equation*}

Hence we can write

\begin{eqnarray*}
\left| R_{q_{\alpha _{k}}^{s}}f\left( x\right) \right| &\geq &\frac{c}{%
\alpha _{k}}\frac{M_{2\alpha _{k}}^{1/p-1}}{\sqrt{\alpha _{k}}}\frac{%
cM_{2s}^{2}}{M_{2\alpha _{k}}}-c\frac{M_{\alpha _{k}}}{\alpha _{k}^{3/2}} \\
&\geq &\frac{cM_{2\alpha _{k}}^{1/p-2}M_{2s}^{2}}{\alpha _{k}^{3/2}}-c\frac{%
M_{\alpha _{k}}}{\alpha _{k}^{3/2}} \\
&\geq &\frac{cM_{2\alpha _{k}}^{1/p-2}M_{2s}^{2}}{\alpha _{k}^{3/2}}.\qquad
\end{eqnarray*}

Then we have

\begin{eqnarray*}
&&\underset{G_{m}}{\int }\left| R^{*}f\left( x\right) \right| ^{p}d\mu
\left( x\right) \\
&\geq &\sum_{s=\left[ 2\alpha _{k}/3\right] }^{\alpha _{k}}\underset{%
I_{2s}\backslash I_{2s+1}}{\int }\left| R_{q_{\alpha _{k}}^{s}}f\left(
x\right) \right| ^{p}d\mu \left( x\right) \\
&\geq &\overset{\alpha _{k}}{\underset{s=\left[ 2\alpha _{k}/3\right] }{\sum
}}\underset{I_{2s}\backslash I_{2s+1}}{\int }\left( \frac{cM_{2\alpha
_{k}}^{1/p-2}M_{2s}^{2}}{\alpha _{k}^{3/2}}\right) ^{p}d\mu \left( x\right)
\\
&\geq &c\sum_{s=\left[ 2\alpha _{k}/3\right] }^{\alpha _{k}}\frac{M_{2\alpha
_{k}}^{1-2p}M_{2s}^{2p-1}}{\alpha _{k}^{3p/2}} \\
&\geq &\left\{
\begin{array}{l}
\frac{2^{\alpha _{k}(1-2p)}}{\alpha _{k}^{3p/2}}\text{,\qquad when }0<p<1/2,
\\
c\alpha _{k}^{1/4},\text{\qquad when }p=1/2,
\end{array}
\right. \rightarrow \infty \text{ ,\qquad when }k\rightarrow \infty .
\end{eqnarray*}

which complete the proof of the theorem 2.

\textbf{Proof of theorem 3.} We write

\begin{eqnarray}
L_{q_{\alpha _{k}}^{s}}f\left( x\right) &=&\frac{1}{l_{q_{\alpha _{k},s}}}%
\underset{j=1}{\overset{q_{\alpha _{k}}^{s}}{\sum }}\frac{S_{j}f\left(
x\right) }{q_{\alpha _{k}}^{s}-j}  \label{16} \\
&=&\frac{1}{l_{q_{\alpha _{k}}^{s}}}\underset{j=1}{\overset{M_{2\alpha
_{k}}-1}{\sum }}\frac{S_{j}f\left( x\right) }{q_{\alpha _{k}}^{s}-j}  \notag
\\
&&+\frac{1}{q_{\alpha _{k}}^{s}}\underset{j=M_{2\alpha _{k}}}{\overset{%
q_{\alpha _{k}}^{s}}{\sum }}\frac{S_{j}f\left( x\right) }{q_{\alpha
_{k}}^{s}-j}  \notag \\
&=&III+IV.  \notag
\end{eqnarray}

Since (see \ref{9})

\begin{equation*}
\left| S_{j}f\left( x\right) \right| \leq c\frac{M_{2\alpha _{k-1}}^{1/p}}{%
\sqrt{\alpha _{k-1}}}\text{ },\text{\qquad }j<M_{2\alpha _{k}}.
\end{equation*}

For $III$ we can write

\begin{equation}
\left| III\right| \leq \frac{c}{\alpha _{k}}\sum_{j=0}^{M_{2\alpha _{k-1}}}%
\frac{1}{q_{\alpha _{k}}^{s}-j}\frac{M_{2\alpha _{k-1}}^{1/p}}{\sqrt{\alpha
_{k-1}}}\text{ }\leq c\frac{M_{2\alpha _{k-1}}^{1/p}}{\sqrt{\alpha _{k-1}}}.
\label{17}
\end{equation}

Using (\ref{11}) we have

\begin{eqnarray}
IV &=&\frac{1}{l_{q_{\alpha _{k}}^{s}}}\sum_{j=M_{2\alpha _{k}}}^{q_{\alpha
_{k}}^{s}}\frac{1}{q_{\alpha _{k},s}-j}\left( \sum_{\eta =0}^{k-1}\frac{%
M_{2\alpha _{\eta }}^{1/p-1}}{\sqrt{\alpha _{\eta }}}\left( D_{M_{_{2\alpha
_{\eta }+1}}}\left( x\right) -D_{M_{_{2\alpha _{\eta }}}}\left( x\right)
\right) \right)  \label{18} \\
&&+\frac{1}{l_{q_{\alpha _{k}}^{s}}}\frac{M_{2\alpha _{k}}^{1/p-1}}{\sqrt{%
\alpha _{k}}}\sum_{j=M_{2\alpha _{k}}}^{q_{\alpha _{k}}^{s}}\frac{\left(
D_{_{j}}\left( x\right) -D_{M_{_{2\alpha _{k}}}}\left( x\right) \right) }{%
q_{\alpha _{k}}^{s}-j}  \notag \\
&=&IV_{1}+IV_{2}.  \notag
\end{eqnarray}

Applying (\ref{4}) in $IV_{1}$ we have
\begin{equation}
\left| IV_{1}\right| \leq c\frac{M_{2\alpha _{k-1}}^{1/p}}{\sqrt{\alpha
_{k-1}}}.  \label{19}
\end{equation}

From (\ref{14}) we obtain

\begin{equation}
IV_{2}=\frac{1}{l_{q_{\alpha _{k},s}}}\frac{M_{2\alpha _{k}}^{1/p-1}}{\sqrt{%
\alpha _{k}}}\psi _{_{M_{2\alpha _{k}}}}\sum_{j=0}^{M_{2s}-1}\frac{%
D_{j}\left( x\right) }{_{M_{2s}}-j}.  \label{20}
\end{equation}

Let $x\in I_{2s}\backslash I_{2s+1}.$ Then $D_{_{j}}\left( x\right) =j,$ $%
j<M_{2s}.$ Consequently

\begin{eqnarray*}
&&.\sum_{j=0}^{M_{2s}-1}\frac{D_{_{j}}\left( x\right) }{M_{2s}-j}%
=\sum_{j=0}^{M_{2s}-1}\frac{j}{M_{2s}-j} \\
&=&\sum_{j=0}^{M_{2s}-1}\left( \frac{M_{2s}}{M_{2s}-j}-1\right) \geq
csM_{2s}.
\end{eqnarray*}

Then

\begin{equation}
\left| IV_{2}\right| \geq c\frac{M_{2\alpha _{k-1}}^{1/p-1}}{\alpha
_{k}^{3/2}}sM_{2s},\qquad x\in I_{2s}\backslash I_{2s+1}.  \label{21}
\end{equation}

Combining (\ref{5}), (\ref{16})-(\ref{21}) for $x\in I_{2s}\backslash
I_{2s+1},s=\left[ 2\alpha _{k}/3\right] ...\alpha _{k}$ and $0<p\leq 1$ we
have

\begin{eqnarray*}
&&\left| L_{q_{\alpha _{k}}^{s}}f\left( x\right) \right| \\
&\geq &c\frac{M_{2\alpha _{k-1}}^{1/p-1}}{\alpha _{k}^{3/2}}sM_{2s}-c\frac{%
M_{\alpha _{k}}}{\alpha _{k}^{3/2}} \\
&\geq &c\frac{M_{2\alpha _{k-1}}^{1/p-1}}{\alpha _{k}^{3/2}}sM_{2s}.
\end{eqnarray*}

Then

\begin{eqnarray*}
&&\int_{G_{m}}\left| L^{*}f\left( x\right) \right| ^{p}d\mu \left( x\right)
\\
&\geq &\sum_{s=\left[ 2\alpha _{k}/3\right] }^{m_{k}}\int_{I_{2s}\backslash
I_{2s+1}}\left| L^{*}f\left( x\right) \right| ^{p}d\mu \left( x\right) \\
&\geq &\sum_{s=\left[ 2\alpha _{k}/3\right] }^{m_{k}}\int_{I_{2s}\backslash
I_{2s+1}}\left| L_{q_{\alpha _{k}}^{s}}f\left( x\right) \right| ^{p}d\mu
\left( x\right) \\
&\geq &c\sum_{s=\left[ 2\alpha _{k}/3\right] }^{m_{k}}\int_{I_{2s}\backslash
I_{2s+1}}\left( \frac{M_{2\alpha _{k-1}}^{1/p-1}}{\alpha _{k}^{3/2}}%
sM_{2s}\right) ^{p}d\mu \left( x\right) \\
&\geq &c\underset{s=\left[ 2\alpha _{k}/3\right] }{\overset{m_{k}}{\sum }}%
\frac{M_{2\alpha _{k-1}}^{1-p}}{\alpha _{k}^{p/2}}M_{2s}^{p-1} \\
&\geq &\left\{
\begin{array}{l}
\frac{2^{\alpha _{k}(1-p)}}{\alpha _{k}^{p/2}}\text{ ,\qquad when }0<p<1, \\
c\sqrt{\alpha _{k}}\text{ ,\qquad when }p=1,
\end{array}
\right. \rightarrow \infty ,\text{ \qquad when }k\rightarrow \infty .
\end{eqnarray*}

Theorem 3 is proved.

\end{document}